\documentstyle[12pt]{article}
\begin{document}
\noindent 
{\large  \bf On the equivalence of two expected average cost
criteria for semi-Markov control processes}\\

\noindent {\bf Anna Ja\'skiewicz} \footnote{Research supported by
KBN Grant 5 P03A 01420}\\
\noindent Institute of Mathematics, Wroc{\l}aw University of
Technology, Wybrze\.ze Wyspia\'nskiego 27, 50-370 Wroc{\l}aw, Poland
(e-mail: ajaskiew@im.pwr.wroc.pl)\\

\noindent {\bf Abstract:} The two expected average costs used in the theory of 
semi-Markov control processes with a Borel state space are
considered. Under some stochastic stability conditions, we prove that the 
two criteria are equivalent in the sense that they lead to the same optimality 
equation.\\

\noindent {\bf Key Words:} Semi-Markov control models,
Borel state space, average cost optimality equation.

\section{The model}
In this paper we study two basic optimality criteria used
in the theory of semi-Markov control processes (see \cite{fein,
minosaki,ross}, for instance).
According to the first one, the average cost is the $\limsup$ of the
expected total costs over a finite number of jumps divided by the
expected cumulative time of these jumps. According to the second
definition, the average cost is the $\limsup$ of the expected total costs 
over the finite deterministic horizon divided by the length of the horizon.
We shall call them (following Feinberg \cite{fein}) the {\it ratio-average 
cost} and {\it time-average cost,} respectively.

Some partial results on the equivalence of the two criteria were given
in the books by Ross \cite{ross} and Puterman \cite{puter}, but only
for countable state space models and stationary policies. A fairly
complete treatment of the problem for semi-Markov control models with
countably many states can be found in \cite{schal}. The main objective
in this paper is to prove the equivalence of the mentioned criteria in
a Borel state space framework. A more detailed presentation of the
perspective in which our research is done is given at the end of this
section.

We shall use the following terminology. A Borel subset, 
say $S,$ of a complete, separable metric space is called a Borel space, 
and it is equipped with the Borel $\sigma$-algebra, denoted by ${\cal B}(S).$ 

A {\it semi-Markov control model} operates as follows. The process is 
observed at time $t=0$ to be in some state $x_0=x\in X,$ where $X$ is 
{\it Borel state space.} At that time an action $a_0=a\in A(x)$ is chosen, 
where $A(x)$ is a compact {\it set of actions available in state} $x.$ The set 
of all actions is $A$ and is also assumed to be a Borel state space. 
By
$$K:= \{(x,a): x\in X, a\in A(x)\}$$  we define the admissible pair set as a 
Borel subset of $X\times A.$

If the current state is $x$ and action $a$ is selected, then the immediate 
cost of $c_1(x,a)$ is incurred and the system remains in state $x_0=x$ for a 
random time $T$ with the cumulative distribution depending only on $x$ and $a.$
 The cost of $c_2(x,a)$ per unit time is incurred until the next transition 
occurs. Afterwards the system jumps to the state $x_1=y$ according to the 
probability measure ({\it transition law}) $q(\cdot|x,a).$ This procedure 
yields a trajectory $(x_0,a_0,t_1,x_1,a_1,t_2,\ldots)$ of some stochastic 
process, where $x_n$ is the state, $a_n$ is control variable and $t_n$ is the 
time of the $n$th transition, $n=0,1,\ldots.$ \\

\noindent
{\sc Remark 1:} In this work we shall slightly abuse the notation. 
Namely, the state and action variables will be denoted by small letters 
$x_n$ and $a_n$ as their values. Other random variables defined on the space 
of all trajectories will be written by means of the capital letters, e.g. 
$T_n$ - the random time of $n$th transition, for $n=1,2,\ldots$ 
with $T_0:=0.$\\

Let $H_{n}$ be the space of {\it admissible histories} up to
the $n$th transition:
$$H_{n} := (K \!\times\! [0,+\infty))^{n} \!\times\! X \quad \mbox{where} 
\quad  H_{0} = X.$$
An element $h_n$ of $H_{n}$ is called a {\it partial history} of the process and is of the form
$$h_{n} :=
(x_{0},a_{0},t_{1},\ldots,x_{n-1},a_{n-1},t_{n},x_{n}).$$

A {\it control policy} (or {\it policy}) is a sequence $\pi = \{\pi_{n}\}$, 
where each $\pi_{n}$ is a conditional probability
$\pi_{n}(\cdot|h_{n})$ on the control set $A(x_{n}),$ given the entire 
history $h_{n}$ such that
$$\pi_{n}(A(x_{n})|h_{n}) = 1 
\quad \forall h_{n} \in H_{n},\quad n=0,1,\ldots. $$
The {\it class of all policies} is denoted by $\Pi.$ 

Let  $F$ be the set of all  Borel measurable 
mappings $f : X \mapsto A$ such that $f(x) \in A(x)$ for each $x
\in X. $ It is well known that $F$ is 
nonempty if the sets $A(x)$ are compact (Corollary 1 in
\cite{bp}). 
A sequence $\pi = \{\pi_{n}\},$ $n = 0,1,\ldots $ is a
{\it (nonrandomized) stationary policy} iff there is some $f
\in F$ 
such that $\pi_n(\cdot|h_{n})$ is concentrated at 
$f(x_{n}) \in A(x_{n})$ for all $h_{n} \in H_{n}$ and $n = 0,1,\ldots.$
Thus any stationary policy $\pi = (f,f,\ldots)$  
can be identified with $f \in F.$ 

Let $(\Omega,\cal F)$ be the measurable space consisting of the sample 
(or trajectory)
space $\Omega := (X \!\times\!A \!\times\! [0,+\infty))^{\infty}$
and the corresponding product $\sigma$-algebra $\cal F$. According to 
the theorem of C. Ionescu Tulcea (Proposition V.1.1 in \cite{neveu} or 
Chapter 7 in \cite{bs}), for each initial state $x_0=x \in X$ and
each policy $\pi \in \Pi$, there exists a unique probability measure 
$P_{x}^{\pi}$ on $\cal F$ such that for all $B \in \cal B$$(A)$, $D \in 
\cal B$$(X)$
and $h_{n} =
(x_{0},a_{0},t_{1},\ldots,x_{n-1},a_{n-1},t_{n},x_{n})$
in $H_{n}$, $n=1,2,\ldots,$
$$P_{x}^{\pi}(x_{0}=x) = 1,$$
$$P_{x}^{\pi}(a_{n} \in B| h_{n}) = \pi_{n}(B|h_{n}),$$
$$P_{x}^{\pi}(x_{n+1} \in D|h_{n},a_{n},t_{n+1}) = q(D|x_{n},a_{n}).$$
Moreover,
$$P_{x}^{\pi}(T_{n+1}-T_{n} \le t|h_{n},a_{n}) = P_{x_n}^{a_n}(T_{n+1}-
T_{n} \le t).$$
We remind that the last equality is a consequence of our assumption that 
the conditional distribution of the difference $T_{n+1}-T_n$ given the
whole history up to the $n$th pair $(x_{n},a_{n})\in K$ depends only on
$x_{n}$ and $a_{n}.$ As usual, by $E^{\pi}_{x}$ we shall denote the
expectation operator with respect to the probability measure $P_{x}^{\pi}.$ 

Further, set 
$\tau(x,a)$  for the {\it mean holding (sojourn) time,} i.e.
$$ \tau (x,a) :=\int_0^\infty tP^a_x(dt)= E^{a}_{x}T.$$
Recall that $T$ is the random time spent in state $x.$

Let $\pi \in \Pi,$ $x \in X$ and $t\ge 0$ be fixed. Put
$$ N(t):=\max \{n\ge 0: T_n\le t\}$$ as the counting process. 
Later on, we shall make some assumptions under which 
$P^\pi_x(N(t)<\infty)=1$ 
(see Remark 2). Now we can define the
two kind of {\it expected average costs} as

\begin{equation}
\label{ar}
J(x,\pi) := \limsup_{n \to \infty} \frac 
{E_{x}^{\pi}\left(\sum_{k=0}^{n-1}
c(x_{k},a_{k})\right)}{E_{x}^{\pi}T_{n}},
\end{equation}
and
\begin{equation}
\label{apt}
j(x,\pi) := \limsup_{n \to \infty} \frac 
{E_{x}^{\pi}\left(\sum_{k=0}^{N(t)}
c(x_{k},a_{k})\right)}{t},
\end{equation}
where
$$ c(x,a):=c_{1}(x,a)+\tau(x,a)c_{2}(x,a).$$

We will need the following assumptions: \\

\noindent   
{\bf B} Basic assumptions: \\
(i) for each $x \in  X$, the set $A(x)$ 
is compact metric space;\\
(ii) for each $x \in X$, 
$c(x,\cdot)$ is lower semicontinuous on $A(x);$\\
(iii) for each $x \in X$ and every Borel set $D \subset X,$
the function $q(D|x,\cdot)$ is continuous on $A(x);$ \\
(iv) for each  $x \in X$, $\tau(x,\cdot)$ is continuous on $A(x),$ and there
exist positive constants $b$ and $B$ such that
$$ 
b \le \tau(x,a) \le B 
$$ 
for all $(x,a) \in K;$\\
(v)  there exist a constant $L > 0$ 
and a Borel measurable function $V:X \mapsto [1, \infty)$ 
such that $|c(x,a)| \le LV(x)$  
for every $(x,a) \in K;$ \\ 
(vi) for each $x \in X,$ the function 
$$ \int_{X} V(y)q(dy|x,\cdot)$$
is continuous on $A(x).$ \\
\noindent
{\bf GE} Geometric ergodicity assumptions: \\
(i) there exists a Borel set $C \subset X$ such that for
some $\lambda \in (0,1)$ and $\eta > 0$, we 
have
$$  \int_{X} \!V(y)q(dy|x,a) \leq \lambda V(x) + \eta 1_{C}(x)$$
for each $(x,a) \in K;$ $V$ is the function introduced in
({\bf B}, v);\\ 
(ii) the function $V$ is bounded on $C$, i.e.
$$ v_{C} := \sup_{x \in C} V(x) < \infty;$$\\
(iii) there exist some $\delta \in (0,1)$ 
and a probability measure $\mu$ concentrated on the Borel set $C$ 
with the property that
   $$   q(D|x,a) \ge \delta\mu (D) $$                         
for each Borel set $D \subset C$, $x \in C$ and $a \in A(x).$\\

For any Borel measurable
function $u:X \mapsto R$ we define the weighted norm as
$$ \|u\|_{V} := \sup_{x \in X}\frac{|u(x)|}{V(x)}.$$ 
By $L_{V}^{\infty}$ we denote the Banach space of all Borel measurable 
functions $u$ for which $\|u\|_{V}$ is finite. 
\\

We also make two additional assumptions on the sojourn time $T$:\\

\noindent
{\bf R} Regularity condition: \\
there exist $\epsilon>0$ and $\beta<1$ such that
$$
P^{a}_{x}(T\le\epsilon)\le \beta $$ for all $x\in C$ and $a\in A(x).$\\
\noindent
{\bf I} Uniform integrability condition: $$
\lim_{t\to\infty}\sup_{x\in C}\sup_{a\in A(x)}P^{a}_{x}(T>t)=0. $$
\\

Assumption ($\bf GE$) is basic for this paper. Inequality 
({\bf GE}, i) is called the "drift inequality" and the set $C$
satisfying ($\bf GE$, iii) is called "small set" \cite{mt1}. 
They imply that the
embedded state process $\{x_{n}\}$ governed by a stationary policy is  
a positive recurrent aperiodic Markov chain with the unique invariant 
probability measure, denoted by $\pi_{f}$ (consult Theorem 11.3.4 and
page 116 in \cite{mt1}). Moreover, $\{x_{n}\}$ is $V$-uniformly ergodic
(Theorem 2.3 in \cite{mt2}), that is, there exist
$\theta > 0$ and  $\alpha \in (0,1)$ such that
$$
\Big|\int_{X} \!u(y)q^{n}(dy|x,f) - \int_{X} \!
u(y)\pi_{f}(dy)\Big| \leq V(x)\|u\|_{V}\theta\alpha^{n}       
$$
for every $u \in L_{V}^{\infty}$  and $x \in X$, $n \geq 1$. Here
$q^{n}(\cdot|x,f)$ denotes the $n$-stage transition probability
induced by $q$ and a stationary policy $f.$ Condition ({\bf GE}) is
often used in the theory of control processes and stochastic games on
Borel state spaces \cite{hl,anea}.

Assumptions ($\bf R$, $\bf I$) are not imposed on the processs, when
we work with the ratio-average criterion (\ref{ar}). However, if we do wish 
do deal with the time-average cost (\ref{apt}), condition ($\bf R$) 
is needed. 
It ensures that the infinite number of transitions does not occur in a finite
interval \cite{ross}. Note also that conditions ($\bf R$), ($\bf I$) 
do not coincide with assumption ($\bf B$, iv). ($\bf R$) implies only that 
$\tau(x,a)>\epsilon(1-\beta)$ for $x\in C.$ The reader who is interested 
in examples is referred to \cite{ross,s}.\\

The literature that deals with semi-Markov control processes under
ratio-average cost is very wide, see e.g. \cite{hlv,minosaki,puter,ross,
schal,s} and
references therein. There are only few papers examining the time-average 
cost \cite{fein,minosaki,puter,ross,schal,yus}. The reason lies in the
fact that it is somewhat easier to study the ratio-average cost.
Generally, these two criteria may have nothing to do with each other.
In other words, they may lead to different cost and optimal policies.
Such situations happen even 
if the state space is countable and there are no ergodic (recurrence) 
properties of the transition probabilities induced by stationary policies. 
Ross \cite{ross} first observed that the two expected costs coincide for 
stationary policies if the embedded Markov chain satisfies some strong 
recurrence condition. Much stronger result was given by Sch\"al \cite{schal}, 
who assumed certain stochastic stability assumptions and proved that the 
optimal expected costs are same in both approaches mentioned above. Moreover, 
Sch\"al showed that the average optimality equation is the same in both cases 
and has a solution. We point out that Sch\"al's paper was devoted to 
semi-Markov 
control processes with countable state space. The optimality equations 
for semi-Markov decision models with ratio-average  criterion (\ref{ar}) and 
Borel state spaces were derived quite recently \cite{hlv,aj}. In this work, 
our goal is to generalize the result of Ross for (uncountable) Borel state 
space. Under the $V$-geometric ergodicity assumption \cite{mt1}, 
we shall prove that criteria (\ref{ar}) and (\ref{apt}) coincide when the 
process is induced by stationary policies. At the same time we show that the 
optimality equation, established in \cite{hlv,aj,va} for the models with the 
cost (\ref{ar}) remains valid for the time-average cost (\ref{apt}).

As in earlier papers \cite{ross,schal}, some parts of our proof 
employs basic facts from renewal theory. Certain consequences of 
$V$-geometric ergodicity given in \cite{mt1,mt2} enable us to apply the 
optional sampling theorem of Doob \cite{nev}, which is the core of the
proof.\\

For convenience of the reader we recall the optimality equation, which is our
point of departure. The proofs are given 
in \cite{hlv,aj,va}. They are based on different methods.\\

\noindent {\sc Proposition:} Let ($\bf B$, $\bf GE$) hold.
Then there exist a function $h\in L^{\infty}_{V}$ and a
constant $g$ such that
\begin{equation}
\label{oe}
h(x)=\min_{a\in A(x)}\left[c(x,a)+\int_{X} h(y)q(dy|x,a)-g\tau(x,a)\right]
\end{equation}
for all $x\in X.$
Moreover, 
$$g=\inf_{\pi\in\Pi}J(x,\pi)=J(x,f^{*}),$$
where $f^{*}\in F$ is a Borel measurable selector of minima on the 
right side of
(\ref{oe}), for each $x\in X.$ \\

\section{Main result}

We begin with presenting our main result in this paper:\\
 
\noindent {\sc Theorem:} Assume ({\bf B},{\bf GE},{\bf R},{\bf I}). Then
\\
(a) \\
$$g=\inf_{\pi\in\Pi} j(x,\pi),$$\\
(b) \\
$$ j(x,f)=J(x,f) \quad\mbox{for each}\quad f\in F.$$ \\

\noindent {\sc Corollary:} Any $f^*\in F$ as in Proposition is average 
optimal with respect to time-average criterion (\ref{apt}).\\

The above results are {\it new.} As already noted, some precedessors of our 
theorem were established in \cite{ross,schal}, but only for semi-Markov
control models with countable state spaces.

Our proof owes much to Ross \cite{ross} and Sch\"al \cite{schal}, especially 
in the parts connected with renewal theory. In order to apply the optional 
sampling theorem, we have to study some consequences of the stochastic 
stability ($V$-geometric ergodicity) assumptions and certain measurability 
issues. \\

\noindent{\bf 1. Some corollaries of "drift inequality."}  
For the set $C$ we define {\it the first return time on C} as 
$$N_{C}:=\min\{ n\ge 1: x_{n}\in C\}.$$
If we do wish to distinguish different return times, we write
$N_C(k)$ for the random time of $k$th visit to $C$: these are defined 
inductively by
$$
N_C(1):=N_C
$$
and
$$
N_C(k):= \min\{ n> N_C(k-1): x_{n}\in C\}.$$
\\

\noindent {\sc Lemma 1:} Let ({\bf GE}, i) hold and let $\{x_{n}\}$ denote
the state space process under arbitrarily fixed policy $\pi \in \Pi.$
Then for each initial state $x \in X$, any function $u \in L^{\infty}_{V}$
and $n \ge 1$ we obtain \\
(a)\\
$$
E_{x}^{\pi}\left(|u(x_{n})|; N_{C}\ge n\right) \le \|u\|_{V}\lambda^{n-1}
\left(\lambda V(x) + \eta 1_{C}(x)\right),
$$
(b)\\
$$
\limsup_{n\to\infty}E_{x}^{\pi} \left(|u(x_{n})|; N_{C}\ge n\right) =0,
$$
(c)\\
$$
E_{x}^{\pi}|u(x_{n})| \le \|u\|_{V}\lambda^{n}V(x) + \eta
\left(1+\ldots+\lambda^{n-1}\right).$$
{\it Proof: } Parts (a) and (c) can be obtained by iteration of 
({\bf GE}, i) on
the set $X\setminus C,$ and $X,$ respectively. Part (b) follows from (a).
$\Box$\\

Our next lemma has a well--known counterpart in the theory of Markov chains
called Comparison Theorem \cite{mt1}. The proof is based on Dynkin's Formula
and proceeds along the same lines as in \cite{mt1}, pages 263-264.\\

\noindent {\sc Lemma 2:} Suppose that nonnegative functions $V,$ $r,$
$s$ satisfy the relationship
$$ \int_{X}V(y)q(dy|x,a)\le V(x)-r(x)+s(x).$$
Then for each $x\in X,$ $\pi\in \Pi$ and any stopping time $\cal S,$ we
have
$$
E^{\pi}_{x}\left[\sum_{k=0}^{{\cal S}-1}r(x_{k})\right]\le V(x)+
E^{\pi}_{x}\left[\sum_{k=0}^{{\cal S}-1}s(x_{k})\right]$$
\\

\noindent {\sc Lemma 3:} If assumption ($\bf GE$, i) holds, then the 
following bounds are satisfied for all $x \in X:$\\
(a) $$E^{\pi}_{x}\left( \sum_{n=0}^{N_{C}-1}V(x_{n})\right)\le \psi(x):=
\frac{1}{1-\lambda}V(x)+ \frac{\eta}{1-\lambda}1_{C}(x),$$\\
(b) $$E^{\pi}_{x} N_{C}\le \phi(x):=
\frac{1}{\ln (1/\lambda)}\left(\ln V(x)+
\frac{\eta}{\lambda}1_{C}(x)\right).$$
{\it Proof:} (a) It follows from Lemma 2 
for the stopping time $N_{C}$ and the functions $r:=(1-\lambda)V,$ 
$s:=\eta 1_{C}.$\\
(b) From Jensen's inequality and ($\bf GE$, i), it follows that
$$  \int_{X} \!\ln V(y)q(dy|x,a) \leq \ln V(x) +\ln\lambda+
\frac\eta\lambda  1_{C}(x)$$ and the rest is obtained by applying
Lemma 2. $\Box$\\

\noindent{\bf 2. An analysis of returns to the set $C.$} 
For any $\pi\in\Pi,$ $t\ge 0$ and $x\in X,$ let
$$
M_{x}^{\pi}(t):=E^{\pi}_{x}\left(\sum_{m=1}^{\infty}
1_{[x_{m}\in C,T_{m}\le t]}\right)=
\sum_{m=1}^{\infty}P^{\pi}_{x}(x_{m}\in C, T_{m}\le t)$$
be the expected number of visits in the small set $C$ 
during the time interval $[0,t].$ 

Define the distribution function $H$ in the following way:
$$
H(t):=\left\{\begin{array}{l@{\quad , \quad}l}
\beta & t\in [0,\epsilon)\\
1 & t\ge \epsilon,
\end{array}\right.
$$
where $\beta$ and $\epsilon$ are taken from assumption ($\bf R$). 
For $t<0$ set $H(t):=0.$
The corresponding renewal function is given by
$$M(t):=\sum_{n=0}^{\infty}H^{n*}(t),$$
where $H^{n*}$ is the $n$-fold convolution of $H$ ($H^{0*}:=1$ on
$[0,+\infty]$ and $H^{0*}:=0$ for $t<0$). By ($\bf R$), we have
$$
P^{\pi}_{x}(T_{N_C}\le t)\le P^{\pi}_{x}(T_{1}\le t)\le H(t),
$$
for each $x\in C.$ Making use of this fact and using 
standard methods from renewal theory, one can show the following facts.\\

\noindent
{\sc Lemma 4:} For $\pi\in\Pi,$ $x\in X$ and $t,h \ge 0,$ we have:\\

(a) $M^{\pi}_{x}(t)\le M(t)<\infty,$ \\

(b) $M^{\pi}_{x}(t+h)-M^{\pi}_{x}(t)\le M(h),$\\

(c) if additionally $z(t)$ is a bounded, nonincreasing, and nonnegative 
function such that 
$\lim_{t\to\infty}z(t)=0,$ then 
$$\lim_{t\to\infty}\frac 1t \int_{0}^{t} z(t-u)M^{\pi}_{x}(du)=0.$$
{\it Proof:} 
(a) This part follows from \cite{ross}, and the fact that
$\{T_{n+1}-T_{n}\},$ $n=0,1,\ldots$ are conditionally independent
random variables given the history of the states and actions process.\\
(b) The proof is similar to that of part (a); see also \cite{black}.\\
(c) Use (a) and (b) and the fact that $\frac{M(t)}{t}\to \frac
1{\epsilon(1-\beta)},$ when $t\to\infty$ (see Key Renewal Theory in 
\cite{ross}). $\Box$\\

\noindent
{\bf 3. The proof of Theorem. } We start with some helpful notation. 
For any $n\ge 1,$ we write 
${\cal F}_{n}$ for the $\sigma$-algebra of all events up to the $n$th state. 
Let $\xi$ be a stopping time relative to $\{{\cal F}_n\}.$ By ${\cal F}_\xi,$ 
we denote the $\sigma$-algebra of all events up to the stopping time $\xi,$ 
i.e.,
${\cal F}_{\xi}:=\{D\in {\cal F}: [\xi=n]\cap D
\in {\cal F}_{n}, \forall n\ge 0 \}.$ 

In this section we accept all our
assumptions. For clear-sighted analysis, we divide the proof into a
sequence of Lemmas.  \\ 

Let $h_{n}= 
(x_{0},a_{0},t_{1},\ldots,a_{n-1},t_{n},x_n)\in H_n.$ 
We put
$h'_{n} :=
(x_{0},a_{0},t_{1},\ldots,a_{n-1},t_{n}).$ For any policy $\pi\in\Pi$ 
and $m\ge 1,$ the conditional policy is formally defined by setting
$$\pi_{n}[h'_{m}](\cdot|h_{n}):=\pi_{n+m}(\cdot|h'_{m},h_{n}).$$
By $E^{\pi [h'_m]}_{x_m},$ we denote the (conditional) expectation operator 
corresponding to the conditional probability measure induced by $\pi [h'_m],$
the transition law $q$ and the holding time distribution. In a similar way, 
we define $\pi[x_0,a_0],$ 
$E^{\pi[x_0,a_0]}_{x_0}$ etc. 

Let $h'_n$ be fixed. Put ${\tilde h}_n:=x_n$ and 
${\tilde h}_m:= (x_n,a_n,t_{n+1},\ldots,t_m,x_m)$ for $m>n.$ 
Identify $(h'_{n},{\tilde h}_m)$ with $h_m.$ With any conditional policy 
$\pi[h'_m]$ we associate the usual policy ${\hat \pi}=\{{\hat \pi}_m\},$ 
where  
$$ {\hat \pi}_m(\cdot|{\tilde h}_m)=\pi_m(\cdot|(h'_n,{\tilde h}_m)) =
\pi_m(\cdot|h_m).$$
Here $x_n$ is treated as the initial state. Below we present a simple 
auxiliary result. \\

\noindent {\sc Lemmma 5:} Let $u:\Omega\mapsto R$ be a Borel 
measurable function such that $E^\pi_x|u|<\infty,$ $x\in X,$ $\pi\in\Pi.$ 
Then \\
(a) $$ u^*(x):=\sup_{\pi\in\Pi}E^{\pi}_{x}u(x_0,a_0,T_1,\ldots)$$
is universally measurable function,\\
(b) 
$$
E^{\pi[h'_n]}_{x_n} u(x_n,a_n,T_{n+1},\ldots)=
E^{\hat \pi}_{x_n} u(x_n,a_n,T_{n+1},\ldots)\le \sup_{\pi\in\Pi}
E^{\pi}_{x_n} u(x_n,a_n,T_{n+1},\ldots).
$$
{\it Proof:} For part (a) see \cite{bs, strauch}. Part (b) is obvious.
$\Box$\\

\noindent
{\sc Lemma 6:} For any $\pi\in\Pi,$ $x\in X$ and $t\ge 0,$ we have:\\
(a) 
$$E^{\pi}_{x} \left( \sum_{n=0}^{N(t)}|c(x_{n},a_{n})|\right)
\le L\psi(x)+L\psi_{C}M(t)<\infty,
$$
 (b) 
$$
E^{\pi}_{x}\left( \sum_{n=0}^{N(t)}\tau (x_{n},a_{n})\right)
\le B\phi(x)+B\phi_{C}M(t)<\infty,
$$ 
with $\psi$ and $\phi$ as defined in Lemma 3 and 
$\psi_{C}:=\sup_{x\in C}\psi(x),$ $\phi_{C}:=\sup_{x\in C}\phi(x)$ 
(see ({\bf GE}, ii)). The constants $L,$ $B$ are from ({\bf B}, iv,v).\\

{\it Proof: } (a)
\begin{eqnarray*}
\lefteqn{E^{\pi}_{x}\left( \sum_{n=0}^{N(t)}|c(x_{n},a_{n})|\right) = 
E^{\pi}_{x}\left(\sum_{n=0}^{N_{C}-1\wedge N(t)}|c(x_{n},a_{n})|
\right)} \\
&&+
\sum_{m=1}^{\infty} E^{\pi}_{x}\left(E^{\pi}_{x} \left[\sum_{n=N_{C}(m)}
^{N_C(m+1)-1\wedge N(t)}|c(x_{n},a_{n})|\Big|{\cal
F}_{N_{C}(m)}\right]\right)=\\
&&E^{\pi}_{x}\left( \sum_{n=0}^{N_{C}-1\wedge N(t)}|c(x_{n},a_{n})|
\right)+\\
&&\sum_{m=1}^{\infty} E^{\pi}_{x}\left(E^{\pi [h'_{m}]}_{x_m} 
\left[\sum_{n=0}^{N_C-1\wedge N(t-T_{m})}|c(x_{n},a_{n})|\right];
x_{m}\in C,T_{m}\le t\right)\le\\
&& L\psi(x)+L\psi_{C}M(t)<\infty,
\end{eqnarray*}
The second equality is due to the strong Markov property generalized
to arbitrary policies (see \cite{schal} for a similar argument). 
The conclusion  follows from Lemmas 3(a), 4(a) and 5(b).
The proof of part (b) is similar to that of part (a). $\Box$\\

\noindent
{\sc Remark 2:} If we replace $|c(x,a)|$ by $|c(x,a)|+1$ in the proof 
of Lemma 6(a), we obtain 
$$E^\pi_x\left(\sum_{n=0}^{N(t)}\left(1+|c(x_{n},a_{n})|\right) 
\right)
\le (L+1)\psi(x)+(L+1)\psi_{C}M(t)<\infty.
$$
Hence, it follows that
$E_x^{\pi}N(t)<\infty$ and consequently $N(t)<\infty$
$P^{\pi}_x$-a.e.\\

\noindent
{\sc Lemma 7:} For any $\pi\in\Pi,$ $x\in X,$ we have:\\
$$
\lim_{t\to\infty} \frac1t E^{\pi}_x V(x_{N(t)+1})=0. $$ 

{\it Proof:} \\
{\it Step 1.} For convenience, we put
$$
w_x^{\pi}(t):=E^{\pi}_x(V(x_{N(t)+1});N_C >N(t)+1),
$$
which can be rewritten as 
\begin{equation}
\label{wzor}
w_x^{\pi}(t)=
E^{\pi}_x\left(\sum_{n=0}^{N_C-2}
\omega(x_n,a_n,t-T_n)\right),
\end{equation}
with
$$
\omega(x,a,t):=E^a_x(V(y);y\not\in C, T>t).
$$
Here $y$ denotes the next state. Recall that $T$ is the sojourn
time in the state $x.$\\
{\it Step 2.} We claim that $w^\pi_x(t)$ is nonincreasing in $t.$
This fact follows immediately from the optional sampling theorem.
For this note that\\
(I)
$ V(x_n)1_{[N_C>n]},$ $n\ge 1$ is supermartingale with respect to 
${\cal F}_{n}$:
\begin{eqnarray*}
E^\pi_x\left(V(x_{n+1})1_{[N_C>n+1]}|{\cal F}_n\right)&\le &
E^\pi_x\left(V(x_{n+1})1_{[N_C\ge n+1]}|{\cal F}_n\right) \\
\le \lambda V(x_{n})1_{[N_{C}\ge n+1]}&\le& V(x_{n})1_{[N_{C}>n]};
\end{eqnarray*}
(II)
$$w^\pi_x(t)\le\sum_{n=1}^{\infty}E^\pi_x\left(V(x_n)1_{[N_C>n]}
\right)\le
\frac{\lambda V(x)+\eta 1_C(x)}{1-\lambda}<\infty;$$
(III)
$$
\lim_{n\to\infty}E^\pi_x\left(V(x_n)1_{[N_C>n]};N(t)\ge n\right)\le
\lim_{n\to\infty}E^\pi_x\left(V(x_n)1_{[N_C>n]}\right)=0.$$
It is easy to see that (II) and (III) follow from Lemma 1(a),(b).

Applying the optional sampling theorem \cite{ws,nev} to the above uniformly 
integrable supermartingale, where $N(t_1)+1\le
N(t_2)+1$ ($t_1<t_2$) are two stopping stopping times, we get
$$w^\pi_x(t_1)\ge w^\pi_x(t_2)$$
for all $x\in X$ and $\pi\in \Pi.$\\
{\it Step 3.}
Put
\begin{equation}
\label{s40}
w_{x}(t):=\sup_{\pi\in\Pi}w^{\pi}_{x}(t)
\end{equation}
for $x\in X,$ $t\ge 0.$ By Lemma 4(b), $w_x(t)$ is universally measurable 
in $x$ for each $t\ge0.$ Note that by (\ref{wzor})
\begin{equation}
\label{s41}
w^\pi_x(t)=E^\pi_x\left(\sum_{n=0}^{N_C-2}\omega(x_n,a_n,t-T_n)\right)=
E^{\pi}_{x}\left(\omega(x,a_{0},t)+1_{[x_1\not\in C, T_{1}\le t]}
w_{x_1}^{\pi[x,a_0]} (t-T_1)\right).
\end{equation}
Denoting $T_1$ by $T$ and $x_1$ by $y,$ we obtain
\begin{equation}
\label{s42}
w^{\pi[x,a_0]}_y(t-T)1_{[y\not\in C,T\le t]}=
w_{y}^{\pi[x,a_0]} (t-T) 1_{[y\not\in C,0\le T\le\frac t2]}+
w_{y}^{\pi[x,a_0]} (t-T) 1_{[y\not\in C, \frac t2<T\le t]}.
\end{equation}
Taking into account (\ref{s40}) and Lemma 5(b), we observe that
$$
w_{y}^{\pi[x,a_0]} (t-T)\le w_y(t-T).$$
This, (\ref{s42}) and the monotonicity of the function 
$t\mapsto w_y(t)$ (Step 2) imply that 
\begin{eqnarray}
\label{s43}
w^{\pi[x,a_0]}_y(t-T)1_{[y\not\in C,T\le t]}&\le &
w_{y} (t-T) 1_{[y\not\in C,0\le T\le\frac t2]}+\\ \nonumber
w_{y} (t-T) 1_{[y\not\in C, \frac t2<T\le t]} &\le& 
w_y\left(\frac t2\right)1_{[y\not\in C]}+w_y(0)1_{[y\not\in C, \frac t2<T]}.
\end{eqnarray}
By Lemma 5(a), both functions $w_y(0)$ and $w_y(\frac t2)$ are universally 
measurable on the state space.
Expressions (\ref{s41}) and (\ref{s43}) yield
\begin{equation}
\label{s44}
w^{\pi}_x(t)\le 
\sup_{a\in A(x)}\omega(x,a,t)+ E^\pi_x \left(w_y\left(\frac t2\right)
1_{[y\not\in C]}\right)+E^\pi_x w_y(0)1_{[y\not\in C, \frac t2<T]}.
\end{equation}
By ($\bf I$),  for $\varepsilon>0$ 
there exists a constant $S>0$ such that for $t>S$
$$
\sup_{x\in C}\sup_{a\in A(x)} P^a_x(T>t)\le \frac{\varepsilon}
{3k[\lambda v_C+\eta]},$$
where $k$ is an integer satisfying 
$$k>\frac 1{\ln\lambda}\ln\frac{\varepsilon(1-\lambda)}
{3[\lambda v_c+\eta]}.$$ 
Let $x\in C$ and $t>2S.$ We obtain the following upper bounds
\begin{equation}
\label{s45}
\sup_{x\in C}\sup_{a\in A(x)} \omega(x,a,t) \le \frac\varepsilon{3k}
\end{equation}
and
\begin{equation}
\label{s46}
E^\pi_x\left(w_y(0)1_{[y\not\in C, \frac t2<T]}\right)\le
\sup_{x\in C}\sup_{a\in A(x)}\left(\int_{X}w_y(0)1_{[y\not\in C]}q(dy|x,a)
\int_{\frac t2}^\infty P^a_x(dt')\right)\le
\frac\varepsilon{3k}.
\end{equation}
From (\ref{s44}), (\ref{s45}), and (\ref{s46}), we get
\begin{equation}
\label{suprema}
\nonumber
w^\pi_x(t)\le 
2 \frac\varepsilon{3k}+\sup_{a\in A(x)}\int_X
w_y\left(\frac t2\right)1_{[y\not\in C]}q(dy|x,a).
\end{equation}
Put 
$$\gamma(y):=w_y\left(\frac t2\right)1_{[y\not\in C]},$$ for $y\in X.$ 
Then  $\|\gamma\|_V<\infty$ and $\gamma$ is universally measurable. 
By Lemma 8.3.7(a) \cite{hl} and F 3.9 \cite{ap}, the function 
$$a\mapsto\int_X \gamma(y)q(dy|x,a)$$ is continuous on $A(x).$ 
Note also, that $$x\mapsto\int_X\gamma(y)q(dy|x,a)$$ is universally
measurable, Proposition 7.46 \cite{bs} or F 3.8 \cite{ap}. Hence, by 
F 2.7 \cite{ap}, there exists a universally 
measurable control function $u$ such that
\begin{equation}
\label{s47}
\max_{a\in A(x)}\int_X\gamma(y)q(dy|x,a)= 
\int_X\gamma(y)q(dy|x,u)
\end{equation}
for each $x\in X.$
By (\ref{suprema}) and (\ref{s47}), we have
$$ w_{x}(t)=\sup_{\pi\in\Pi}w^\pi_x(t)\le 
2 \frac\varepsilon{3k}+\int_X
w_y\left(\frac t2\right)1_{[y\not\in C]}q(dy|x,u).
$$
Iteration of the last inequality $(k-1)$ times together with Step 2(II)
and
Lemma 1(a) (which is also valid for universally measurable policies) 
gives
\begin{eqnarray}
\label{s48}
\sup_{x\in C}\sup_{\pi\in \Pi}w^\pi_x(t)&\le& 
\frac{2\varepsilon}{3}
+\sup_{x\in C}E^{u}_{x}\left(w_{x_k}\left(\frac
t{2^k}\right);N_C>k\right)\\ \nonumber
\le \frac{2\varepsilon}{3}
&+&\sup_{x\in C}E^{u}_{x}\left(\frac{\lambda}{1-\lambda} 
V(x_k);N_C>k\right)\le \varepsilon.
\end{eqnarray}
{\it Step 4.} Set
$$
z_x^{\pi}(t):=E^{\pi}_x(V(x_{N(t)+1});N_C\ge N(t)+1),
$$
and
$$
z(t):=\sup_{x\in C}\sup_{a\in A(x)}z^{\pi}_x(t).
$$
Proceeding analogously as in Step 2, we note that $z(t)$ is a nonincreasing 
function in $t$ and $z(t)<\infty.$ Moreover, using the strong Markov property 
and Lemma 5(b), we obtain
\begin{eqnarray*}
\lefteqn
{E^{\pi}_x V(x_{N(t)+1})= z^{\pi}_x(t)+}\\
&&\sum_{m=1}^{\infty} E^{\pi}_x\left(
E^{\pi}_x \left[V(x_{N(t)+1});N_C(m+1)\ge N(t)+1>N_C(m)\Big|{\cal F}
_{N_C(m)}\right]\right) \\
&&\le 
z^{\pi}_x(t)+\\
&&\sum_{m=1}^{\infty} E^{\pi}_x\left(
E^{\pi[h'_{m}]}_{x_m}\left[V(x_{N(t-T_m)+1});
N_C\ge N(t-T_m)+1\right];x_m\in C, T_m\le t\right)\\
&&\le 
z^{\pi}_x(t)+\\
&&\sum_{m=1}^{\infty} E^{\pi}_x\left(
\sup_{{\bar x}\in C}\sup_{{\hat \pi}\in\Pi}E^{\hat\pi}_{{\bar x}} 
\left[V(x_{N(t-T_m)+1});
N_C\ge N(t-T_m)+1\right];x_m\in C, T_m\le t\right)\\
&&=
z^{\pi}_x(t)+\int_0^t z(t-u)M^{\pi}_x(du).
\end{eqnarray*}
{\it Step 5.} Now it remains only to prove that $\lim_{t\to\infty}z(t)=0.$
This follows from
$$
z(t)\le \sup_{x\in C}\sup_{\pi\in \Pi}w^\pi_x(t)+\sup_{x\in C}
\sup_{\pi\in \Pi}E^\pi_x\left(V(x_{N(t)+1});N_C=N(t)+1\right),
$$
and
\begin{eqnarray*}
\sup_{x\in C}\sup_{\pi\in
\Pi}E^\pi_x\left(V(x_{N_C});N_C=N(t)+1\right)
&\le& v_C \sup_{x\in C}
\sup_{\pi\in \Pi}P^\pi_x(T_{N_C}>t)\\
\le v_C \sup_{x\in C}
\sup_{\pi\in \Pi}\frac{E^\pi_x(T_{N_C})}{t}&\le& \frac{B\phi_C}{t}.
\end{eqnarray*}
The sequence of inequalities is due to assumption ({\bf GE},ii), Markov 
inequality and Lemma 1(b) (see also ($\bf B$, iv)), respectively. Hence, 
by (\ref{s48}) and the last expression 
$z(t)\to 0$ as $t\to\infty.$ 

The desired assertion is a consequence of
$$E^{\pi}_x V(x_{N(t)+1})\le z^{\pi}_x(t)+\int_0^t z(t-u)M^{\pi}_x(du) 
\quad \mbox{(by Step 4),}$$
Lemma 4(c) and the fact that 
$$z^\pi_x(t) \le \frac{\lambda+\eta}{1-\lambda}V(x),$$
which follows from Lemma 1(a).  $\Box$\\

\noindent
{\sc Lemma 8:} For any $\pi\in\Pi,$ $x\in X,$  we have:\\

$$\lim_{t\to\infty} \frac 1t E^{\pi}_x T_{N(t)+1}=1.$$\\

{\it Proof:} Note that 
$$
\frac{t}{t}\le
E^{\pi}_x \frac{T_{N(t)+1}}{t}\le \frac{E^{\pi}_x 
T_{N(t)}+E^{\pi}_{x}(T_{N(t)+1}-T_{N(t)})}{t}\le \frac{t+B}{t},$$
where $B$ is from assumption ($\bf B$, iv).
The last inequality is due to the fact that $[N(t)+1=n]\in {\cal
F}_{n}$ and
\begin{eqnarray*}
E^{\pi}_{x}(T_{N(t)+1}-T_{N(t)})&=&\sum_{n=1}^{\infty}
E^{\pi}_{x}\left(T_{n+1}-T_{n}\right)1_{[N(t)+1=n]}\\ &=&
\sum_{n=1}^{\infty}
E^{\pi}_{x}\left(E^{\pi}_{x}\left[\left(T_{n+1}-T_{n}\right)
1_{[N(t)+1=n]}\Big|
{\cal F}_{n}\right]\right)\\
&=&
\sum_{n=1}^{\infty}
E^{\pi}_{x}\left(\tau(x_{n},a_{n})1_{[N(t)+1=n]}\right)\le
B. \quad \Box
\end{eqnarray*}
\\

{\it Proof of Theorem:} 
(a) We claim that 
$$S_{n}:=\sum_{k=0}^{n-1}\left(c(x_{k},a_{k})-g\tau(x_{k},a_{k})\right)
+h(x_{n})$$
is a submartingale with respect to ${\cal F}_{n}.$  This follows from
the optimality equation (\ref{oe}), because 
$$ h(x)\le c(x,a)-g\tau(x,a)+\int h(y)q(dy|x,a).$$
In order to apply optional sampling theorem for this submartingale,
where $N(t)+1$ is a stopping time, we have to check that\\
(I) $E^{\pi}_{x}|S_{N(t)+1}|$ is well defined; \\
(II) $E^{\pi}_{x}[|S_{n}|; N(t)\ge n]$ tends to zero, when
$n\to\infty.$\\
For (I), it holds
$$E^{\pi}_{x}|S_{N(t)+1}|\le
E^{\pi}_{x}\left(\sum_{k=0}^{N(t)}|c(x_{k},a_{k})|\right)+gE^{\pi}_{x}
\left(\sum_{k=0}^{N(t)}\tau(x_{k},a_{k})\right)+ 
\|h\|_{V}E^{\pi}_{x}V(x_{N(t)+1}).
$$
These expressions are finite by Lemmas 6 and 7 ($h\in L_{V}^{\infty}$).
Furthermore,
\begin{eqnarray*}
\lefteqn{E^{\pi}_{x}\left(|S_{n}|; N(t)\ge n\right)\le }\\
&&E^{\pi}_{x}\left(\sum_{k=0}^{N(t)}\left(|c(x_{k},a_{k})|+g
\tau(x_{k},a_{k}); N(t)\ge n\right)\right)+ 
\|h\|_{V}E^{\pi}_{x}(V(x_{n}); N(t)\ge n).
\end{eqnarray*}
Taking into account Remark 2, these terms go to zero by Lemma 6 
(the first one) and by Lemma 1(b) (the second one). Finally, we obtain
$$
h(x)\le E^{\pi}_{x}\left(\sum_{k=0}^{N(t)}\left(c(x_{k},a_{k})-
g\tau(x_{k},a_{k})\right)\right)
+E^{\pi}_{x}h(x_{N(t)+1})$$
and consequently,
$$
g \frac 1t E^{\pi}_{x} T_{N(t)+1}\le  \frac 1t E^{\pi}_{x}\left(
\sum_{k=0}^{N(t)}c(x_{k},a_{k})\right) +E^{\pi}_{x}\left(\frac 
{h(x_{N(t)+1})}{t}\right) -\frac
{h(x)}{t}.
$$
The left side tends to $g$ (Lemma 8), whilst the right side goes to
$j(x,\pi),$ defined in (\ref{apt}) (Lemma 7). \\
(b) Let $f\in F$ be fixed. Then there exists a function $h_{f}\in
L^{\infty}_{V}$ for which Poisson's equation  holds, i.e.
$$
h_{f}(x)=c(x,f)+\int_{X} h_{f}(y)q(dy|x,f)-J(x,f)\tau(x,f).
$$
The arguments used above in particular imply that 
$$S'_{n}:=\sum_{k=0}^{n-1}\left(c(x_{k},f)-J(x,f)\tau(x_{k},f)\right)
+h_{f}(x_{n})$$
is a uniformly integrable martingale. (We recall that by our assumption
($\bf GE$), the ratio-average cost is independent of the initial state for 
each stationary policy \cite{hl}.) Applying Doob's theorem, we get 
$$h_{f}(x)= E^{f}_{x}\left(\sum_{k=0}^{N(t)}\left(c(x_{k},f)-
J(x,f)\tau(x_{k},f)\right)\right)
+E^{f}_{x}h_{f}(x_{N(t)+1}).$$
This gives the result.
$\Box$\\

\noindent
{\sc Remark 3:} In the proof it is assumed that $\tau(x,a)<B$ ($\bf B$, iv). 
However, the optimality equation (\ref{oe}) remains true, if we allow for 
unbounded mean holding time, i.e.
$$
\tau(x,a)<B_{1}V(x),
$$
for some constant $B_{1}.$ The direct proof is provided in
\cite{va}. The reader can also follow the proof given in \cite{aj} with
slight modification of the constants in Theorem 1. Then, the minor
corrections in the proof of Lemmas 5(b) and 7 give the equivalence of 
expected average costs, (\ref{ar}) and (\ref{apt}), 
for stationary policies.\\

\noindent
{\sc Remark 4:} Our main theorem has some relevance to studying stochastic 
games with Borel state space. Namely, the results given in \cite{siam,appl} 
for semi-Markov games remains also valid for time-average criterion 
(\ref{apt}).\\

\noindent
{\bf Acknowledgment:} I wish to thank Prof. Eugene Feinberg and 
Prof. Andrzej S. Nowak for suggesting the problem, many helpful discussions 
and strong encouragment.

\small
\bibliographystyle{abbrv}

\begin{thebibliography}{99}

\bibitem{bs} {\sc D.P. Bertsekas and S.E. Shreve,} {\it Stochastic Optimal 
Control: The Discrete Time Case,} Academic Press, New York, 1978.

\bibitem{black} {\sc D. Blackwell,} {\it A renewal theorem,} Duke Math.
J., 15 (1948), pp. 145--150.

\bibitem{bp} {\sc L.D. Brown and R. Purves,} {\it Measurable selections
of extrema,} Ann. Stat., 1 (1973), pp. 902--912.

\bibitem{fein} {\sc E.A. Feinberg,} {\it Constrained semi-Markov decision
processes with average rewards,} Math. Methods Oper. Res., 39 (1994), 
pp. 257--288.
     
\bibitem{hl} {\sc O. Hern\'{a}ndez-Lerma  and J.B. Lasserre,} {\it Further
     Topics on Discrete-Time Markov Control Process,} Springer-Verlag,
     New York, 1999.

\bibitem{hlv} {\sc O. Hern\'{a}ndez-Lerma and F. Luque-V\'{a}squez,} 
     {\it Semi-Markov control models with average costs,} Applicationes
     Mathematicae, 26 (1999), pp. 315-331.

\bibitem{aj} {\sc A. Ja\'skiewicz,} {\it An approximation approach to
ergodic semi-Markov control processes,} Math. Methods  Oper. Res., 54
(2001), pp. 1--19.

\bibitem{siam} {\sc A. Ja\'skiewicz,} {\it Zero-sum semi-Markov games,} 
SIAM J. Control Optim., 41 (2002), pp. 723--739.

\bibitem{ws} {\sc A.P. Maitra and W.D.Sudderth,} {\it Discrete Gambling and
Stochastic Games,} Springer-Verlag, New York,1996.

\bibitem{mt1} {\sc S.P. Meyn and R.L. Tweedie,} {\it Markov Chains and
Stochastic Stability,} Springer-Verlag, New York, 1993.

\bibitem{mt2} {\sc S.P. Meyn and R.L. Tweedie,} {\it Computable 
bounds for geometric convergence rates of Markov chains,} Ann. Appl. 
Probab., 4 (1994), pp. 981-1011.

\bibitem{minosaki} {\sc H. Mine and S. Osaki,}  {\it Markovian Decision
Processes,} Elsevier, New York, 1970.

\bibitem{neveu} {\sc J. Neveu,} {\it Mathematical Foundations of the Calculus
     of Probability,} Holden-Day, San Francisco, 1965.

\bibitem{nev} {\sc J. Neveu,} {\it Discrete-Parameter Martingales,}
     Elsevier, New York, 1975.

\bibitem{ap} {\sc A.S. Nowak,} {\it Universally measurable strategies 
in zero-sum stochastic games,} Ann. Probab., 13 (1985), pp. 269--287.

\bibitem{appl} {\sc A.S. Nowak,} {\it Some remarks on equilibria in 
semi-Markov games,} Applicationes Mathematicae, 27 (2000), pp. 385--394.

\bibitem{anea} {\sc A.S. Nowak and E. Altman} {\it
$\varepsilon$-Equilibria for stochastic games with uncountable state
space and unbounded costs,} SIAM J. Control Optim., 40 (2002), 
pp. 1821--1839.


\bibitem{puter} {\sc L.M. Puterman,} {\it Markov Decision Processes,}
John Wiley, New York, 1994. 

\bibitem{ross} {\sc S.M. Ross,} {\it Applied Probability Models with
     Optimization Applications,} Holden-Day, San Francisco, 1970.

\bibitem{schal} {\sc M. Sch\"al,} {\it On the second optimality equation for 
semi-Markov decision models,} Math. Oper. Res., 17 (1992), pp. 470--486.

\bibitem{s} {\sc L.I. Sennott,}  {\it Average cost semi-Markov decision
     processes and the control of queueing system,} Probability in the
     Engineering and Informational Sciences, 3 (1989), pp. 247--272.

\bibitem{strauch} {\sc R.E. Strauch,} {\it Negative dynamic
programming,} Ann. Stat., 37 (1966), pp. 871--890.

\bibitem{yus} {\sc A. Yushkevich,} {\it On semi-Markov controlled models
with an average reward criterion,}  Theory Probab. Appl., 26 (1981), 
pp. 796--803.

\bibitem{va} {\sc O. Vega-Amaya and F. Luque-V\'{a}squez,} 
{\it Sample-path average cost optimality for semi-Markov control 
processes on Borel spaces: unbounded costs and mean holding times,} 
Applicationes Mathematicae, 27 (2000), pp. 343--367.

\end{thebibliography}

\end{document}